\documentclass[12pt]{report}
\usepackage{type1ec}
\usepackage[T2A]{fontenc}
\usepackage[russian]{babel}

\usepackage{amsmath,amsfonts,amssymb,latexsym,amsthm,fancyhdr,multicol}
\usepackage[a4paper,pdftex,dvips,nofoot,verbose,
            headheight=17pt,headsep=0pt,
            textheight=250mm,bottom=20mm,
            left=20mm,right=20mm
            ]{geometry}

\usepackage[pdftex]{graphicx}
\usepackage{url}

 \makeatletter
\renewcommand{\@oddhead}{\hfill    ~\thepage~   \hfill\vspace{0.5 cm}}

\newcounter{tabular}

\begin{document}
\large

\begin{center}\Large\bf
О свойствах коэффициентов экзотических и сложных разложений решений шестого уравнений Пенлеве
\end{center}
\begin{center}
\large\bf И.\,В.\,Горючкина
\end{center}

\par\medskip
{\normalsize Известно, что среди формальных решений шестого уравнения Пенлеве встречаются ряды по целым степеням независимой комплексной переменной $x$
 с коэффициентами в виде формальных рядов Лорана (с конечными главными частями) от $\ln^{-1}{x}$ \,({\it сложные разложения})\,
 или  $x^{\,{\rm i}\, \theta}$, где ${\rm i}=\sqrt{-1},\;\theta\in\mathbb{R},$ $\theta\neq 0$ \,({\it экзотические разложения}).
 Эти коэффициенты можно вычислять последовательно. Здесь исследуются аналитические свойства рядов, которые являются коэффициентами
     сложных и экзотических разложений решений шестого уравнения Пенлеве.}

\begin{center}
{\Large\bf On properties of the coefficients of the complicated and exotic expansions of the solutions of the sixth Painlev\'e equation}
\end{center}
\begin{center}
 \bf I.\,V.\,Goryuchkina
\end{center}

\par\medskip
{\normalsize It is known, that among the formal solutions of the sixth Painlev\'e equation there met series with integer power exponents of the independent variable $x$
with coefficients in form of formal Laurent series (with finite main parts) in $\log^{-1}{x}$ \,({\it complicated expansions}),\,    or in $x^{\,{\rm i}\, \theta}$,
where ${\rm i}=\sqrt{-1},\;\theta\in\mathbb{R}$, $\theta\neq 0$ \,({\it exotic expansions}). These coefficients can be computed consecutively. Here we research analytic properties of the series,
that are the coefficients of the complicated and exotic expansions of the solutions of the sixth Painlev\'e equation.}

\bigskip
{\bf 1. Постановка задачи.} Рассматривается шестое уравнение Пенлеве
$$y''=\frac{(y')^2}{2}\left(\frac{1}{y}+\frac{1}{y-1}+\frac{1}{y-x}\right)-y'\left(\frac{1}{x}
+\frac{1}{x-1}+\frac{1}{y-x}\right)+$$
\begin{equation}+\frac{y(y-1)(y-x)}{x^2(x-1)^2}\left[a+b\frac{x}{y^2}+c\frac{x-1}{(y-1)^2}+
d\frac{x(x-1)}{(y-x)^2}\right],\label{eq1}\end{equation}где  $a,\;b,\;c,\;d$ --
комплексные параметры, $\,x\,$  и  $\,y\,$ -- комплексные переменные,
$y'=dy/dx$. Оно имеет три особые точки независимой переменной $\;x=0,\;x=1\,$  и  $\,x=\infty$,
и симметрии уравнения, порожденные тремя преобразованиями $$1)\, x=z,\;y=z/w,\quad 2)\, x=1/z,\;y=1/w,\quad 3)\, x=1-z,\;y=1-w,$$ которые позволяют
переводить формальные решения вблизи одной особой точки в формальные решения вблизи той же самой или другой особой точки уравнения (см. [\ref{Gromak}]). Эти симметрии сильно упрощают вычисления. Формальные решения в окрестности точки $x=0$ уравнения \eqref{eq1}, из которых с помощью симметрий уравнения получаются все остальные его формальные решения в окрестности всех трех его особых точек,  называются {\it базовыми} формальными решениями шестого уравнения Пенлеве.

Шестое уравнение Пенлеве имеет формальные решения в виде степенных рядов независимой переменной
 $\,x\,$  с коэффициентами, которые являются многочленами Лорана или формальными рядами Лорана с конечными главными частями от
  $\;\ln^{-1} x\;$ ({\it степенно-логарифмические \,{\rm или}\, сложные разложения}) и $\;x^{\,{\rm i}\,\theta}\;$ ($\theta\in\mathbb{R}$,
  ${\rm i}=\sqrt{-1}$)\, ({\it степенные\, {\rm или}\, экзотические разложения}) с комплексными коэффициентами,  см. [\ref{MMO}]. Далее мы будем изучать только формальные решения шестого уравнения Пенлеве,
  которые имеют вид рядов по целым степеням $\,x\,$ с\, коэффициентами, которые являются рядами от $\;\ln^{-1} x\;$ или $\;x^{\,{\rm i}\,\theta}\;$ ($\theta\neq 0$).\,  Но, как уже было сказано, в силу симметрий уравнения  \eqref{eq1} изучение свойств формальных решений этого уравнения сводится к изучению свойств только его базовых формальных решений, а свойства всех остальных его формальных решений будут подобными. Все интересующие нас базовые формальные решения имеют вид
\begin{equation}\label{eq2}
  y=\sum\limits_{k=0}^{\infty}\varphi_k(x)\,x^k,
\end{equation}
где $\varphi_k(x)$ -- это формальные ряды вида
$$
\varphi_k(x)=\xi^\alpha\,\sum\limits_{j=0}^{\infty}\,c_{kj}\,\xi^j, \;\;\qquad c_{kj}\in\mathbb{C}, \;\;\alpha\in\mathbb{Z},
$$
\begin{equation}\label{eq9}
\xi=\ln^{-1}{x}\qquad \mbox{ или }\qquad\xi=x^{\,{\rm i}\,\theta}\;\;(\theta\neq 0).
\end{equation}
\medskip
Сначала мы будем проводить анализ формальных решений шестого уравнения Пенлеве, используя общее выражение \eqref{eq2}, а затем -- уже конкретные
формальные решения шестого уравнения Пенлеве, вычисленные ранее с помощью методов плоской степенной геометрии [\ref{MMO}].

 Отметим, что в работе [\ref{OnClasses}] построены общие замкнутые относительно алгебраических операций и дифференцирования классы формальных
 решений конечного порядка алгебраических (полиномиальных) обыкновенных дифференциальных уравнений, которые могут быть вычислены с помощью
  методов плоской степенной геометрии. Кроме того, в [\ref{OnClasses}] доказывается теорема о том, что {\it если существует (выделенного класса)
  формальное решение алгебраического обыкновенного дифференциального уравнения, то первое приближение этого решения является решением первого
  приближения этого уравнения}.\, При этом первые приближения решения и уравнения удобно выделять посредством многоугольника Ньютона-Брюно,
  построение которого основано на понятиях порядка функции и порядка формального ряда. Результаты работы [\ref{OnClasses}] будут использоваться
  в этой работе. А именно, мы будем пользоваться замкнутостью выделенных в работе [\ref{OnClasses}] классов формальных решений алгебраического
  обыкновенного дифференциального уравнения. Стоит сказать, что к выделенному классу формальных решений относятся все формальные решения
  шестого уравнения Пенлеве.
   В этой работе доказывается теорема (подобная теореме из [\ref{OnClasses}], формулировка которой выделенна курсивом выше)
   для формальных решений того же класса, что и в работе [\ref{OnClasses}], но уже для обыкновенных дифференциальных уравнений
   полиномиальных по зависимой переменной и ее производным, но с коэффициентами в виде рядов конечного порядка. Например, уравнение такого класса получается
     после замены зависимой переменной $\;y=\varphi_0(x)+u\;$ в алгебраическом обыкновенном дифференциальном уравнении,
     где $\;\varphi_0(x)\;$ --\, это ряд конечного порядка. Применяя эту теорему к новому (не алгебраическому обыкновенному дифференциальному) уравнению,
      можно выделить из
      него  приближенное уравнение,
    решением которого является ряд $\;u=\varphi_1(x) x$.\, Таким образом, эта теорема используется при вычислении дальнейших членов ряда \eqref{eq2},
    который является решением алгебраического обыкновенного дифференциального уравнения. Понятия порядка функции и
   порядка формального ряда дано в следующем разделе.

Также отметим, что мы не можем предъявить общие или рекурсивные формулы для коэф\-фи\-циен\-тов $\varphi_k(x)$ формальных решений вида \eqref{eq2}
уравнения \eqref{eq1}, но тем не менее, мы можем исследовать уравнения на эти коэффициеты, которые линейны.
А значит, пользуясь теоремами теории аналитических дифференциальных уравнений, можем судить о свойствах решений этих уравнений.
Тоже самое имеет место для коэффициентов формальных решений других алгебраических обыкновенных дифференциальных уравнений.
Уравнения на коэффициенты нелинейны лишь когда первая вариация приближенного уравнения вдоль его решения равна нулю (т.~е. когда приближенное
уравнение имеет кратное решение), что скорее является исключением.
Поэтому метод исследования достаточно общий.

Здесь и далее используется одно и тоже обозначение $\varphi_k(x)$ для формального ряда и его суммы.

Шестое уравнение Пенлеве имеет двупараметрическое (неформальное) решение с асимптотиками $\varphi_0(x)$ в виде рациональных функций
 от $x^{\,{\rm i}\,\theta},$\, $\theta\neq 0$ в окрестности всех особых точек ($x=0,\;1,\;\infty$) уравнения.
Согласно работе [\ref{Guzzetti}] нули и полюсы настоящего решения сосредоточены в областях, в которых лежат нули и полюсы этих асимптотик
(т.~е. в зеркальных секторах с вершинами в особых точках с биссектрисами на которых накапливаются полюса этих асимптотик).
Результаты этой и проделанной автором работы согласуются, т.е. коэффициенты $\varphi_k(x)$ экзотических разложений имеют особые точки в этих же секторах.

Ранее автором была доказана теорема о сходимости формальных решений алгебраических ОДУ (см. [\ref{OnConv}], [\ref{GG}])
согласно которой все  формальные решения шестого уравнения Пенлеве, имеющие вид степенных рядов с постоянными комплексными
коэффициентами и комплексными показателями степени (вещественные части которых частично упорядочены,  не имеют точек накопления и
при фиксированном значении вещественной части имеется конечное число
показателей степени с одинаковой вещественной частью),
 сходятся в некоторых областях комплексной плоскости. Эту теорему можно также применять к коэффициентам сложных и экзотических разложений
 решений шестого уравнения Пенлеве.

В этой работе решается следующая задача: {\it исследовать асимптотическое поведение и аналитические свойства коэффициентов $\varphi_k(x)$
формальных решений вида \eqref{eq2} уравнения \eqref{eq1}.}

Но перед тем как приступить к изложению основных результатов, сформулируем и докажем теорему о первом приближении
 формального решения обыкновенного дифференциального уравнения, которая, как уже было ранее написано, отличается от подобной теоремы из [\ref{OnClasses}]
 классом уравнения. А именно: в работе [\ref{OnClasses}] рассматривалось алгебраическое  обыкновенное дифференциальное уравнение,
 а здесь мы будем рассматривать обыкновенное дифференциальное уравнение полиномиальное относительно зависимой переменой и ее производных,
 конечного набора степенных функций с комплексными степенями и конечного набора формальных рядов Лорана с конечной главной частью от функций нулевого порядка.

\bigskip
{\bf 2. Теорема о первом приближении.} Поскольку понятия порядка функции и порядка ряда  нам необходимы в дальнейшем и уже употреблялись в предыдущем параграфе, напомним эти понятия.

 {\it Функция} $\;\varphi(x)\;$ в точке $\;x=0\;$  {\it имеет порядок  $\;r\in \mathbb{R}\bigcup\{\infty\}$,}\, если существует предел  \begin{equation}\lim\limits_{\scriptstyle x\rightarrow 0\atop \scriptstyle x\in\mathcal{D}}\frac{\ln{|\varphi(x)|}}{\displaystyle\ln{|{x}|}}=r,\label{eq2n}\end{equation}
 где $\;\mathcal{D}\;$ --  это  открытое связное подмножество множества комплексных чисел, замыкание которого содержит нуль.
 {\it Функция} $\;\varphi(x)\;$ в точке $\;x=\infty\;$  {\it имеет порядок  $\;r\in \mathbb{R}\bigcup\{\infty\}$,}\, если существует предел  \begin{equation}\lim\limits_{\scriptstyle x\rightarrow \infty\atop \scriptstyle \;x\in\mathcal{D}_{\infty}}\frac{\ln{|\varphi(x)|}}{\displaystyle\ln{|{x}|}}=r,\label{eq2nn}\end{equation}
 где  $\;\mathcal{D}_{\infty}\;$ -- это открытое связное подмножество множества комплексных чисел, замыкание которого содержит бесконечность.

{\it Функция $\;\varphi(x),\;$ имеющая порядок $\;r\in\mathbb{R}\bigcup\{\infty\}\;$ и в нуле, и в бесконечности, называется функцией  порядка $\,r$.}

{\it Порядком формального ряда вида \eqref{eq2} будем считать
порядок его первого ненулевого члена.}

Рассмотрим уравнение
\begin{equation}\label{eq25}
g(x,u,u',\dots,u^{(n)})=0,
\end{equation}
 левая часть которого является суммой членов вида
\begin{equation}\label{eq26}
\psi(x)\,x^{q_1}u^{q_{20}}(xu')^{q_{21}}\dots(x^nu^{(n)})^{q_{2n}},
\end{equation}
где $\psi(x)$ -- функция нулевого порядка или формальный ряд Лорана с конечной главной частью от функций нулевого порядка (т.е. это ряд многих переменных), $q_1\in\mathbb{C}$, $q_{21},\dots,q_{2n}\in\mathbb{Z}_+$, и если $\psi^{(n)}(x)\not\equiv 0$, то порядок формального ряда $\psi^{(n)}(x)$ равен $-n$.

Пусть уравнение \eqref{eq25} имеет формальное решение $u=\varphi$ в виде ряда
\begin{equation}\label{eq27}
\varphi=\sum\limits_{k=0}^{\infty}\varphi_k(x)x^{\rho_k},
\end{equation}
где $\varphi_k(x)$ -- функции нулевого порядка или ряды Лорана с конечными главными частями от функции нулевого порядка,  $\rho_k\in\mathbb{R}$, и если $\varphi_k^{(n)}(x)\not\equiv 0$, то порядок $\varphi_k^{(n)}(x)$ равен $-n$, при этом, $\rho_{k+1}>\rho_k$.

Каждому слагаемому вида \eqref{eq26} уравнения \eqref{eq25} можно поставить в соответствие точку (векторный показатель степени этого слагаемого) $\,(q_1,\,q_2)$,\, $q_2=q_{20}+{\dots}+q_{2n}.$ Согласно результатам работы [\ref{OnClasses}] порядок ряда $$\displaystyle \psi(x)\,x^{q_1}\varphi^{q_{20}}(x\varphi')^{q_{21}}\dots(x^n\varphi^{(n)})^{q_{2n}}$$ равен $q_1+q_2\rho_0$.

Пусть множество $\,\{Q_i=(q_1^i,\,q_2^i),\;i=0,\dots,m\}$\, -- это носитель уравнения \eqref{eq25}, т.~е. множество, содержащее все векторные показатели степени всех слагаемых этого уравнения, а $\;R\;$ -- это вектор $\;(1,\rho_0)$.
Рассмотрим все возможные скалярные произведения $\;\langle Q_i,\,R \rangle=c_i\in\mathbb{R}.$\, Пусть $\;c=\min\limits_{i=0,\dots,m}\, c_i$.\, Сумма слагаемых вида \eqref{eq26} уравнения \eqref{eq25}  с такими векторными показателями степени, что $\;\langle Q_i,\,R \rangle=c\in\mathbb{R},\;$ называется {\it укороченной суммой} (см.  [\ref{BrunoUMN}]), обозначим ее $\;\hat{g}(x,u_0,\dots,u_n)$,\, а уравнение
\begin{equation}\label{eq28}
\hat{g}(x,u_0,\dots,u_n)=0
\end{equation}
 -- {\it укороченным уравнением}.

\bigskip\bigskip
{\bf Теорема 1. } \,{\it Если уравнение \eqref{eq25} имеет формальное решение $u=\varphi$, где $\varphi$ определено формулой \eqref{eq27}, тогда укороченное уравнение \eqref{eq28} имеет укороченное $($фор\-маль\-ное$)$ решение
\begin{equation}\label{eq44}
u=\hat{\varphi},\quad\quad\hat{\varphi}=\varphi_0(x).
\end{equation}
}

\medskip
{\bf Доказательство.} \,Если $\;\rho_0\neq 0,\;$ то в уравнении \eqref{eq25} и его формальном решении ${\varphi}$  сделаем степенное преобразование
\begin{equation}u=x^{\rho_0}v.\label{eq29}\end{equation} Получаем уравнение \begin{equation} G(x,v,v',\dots,v^{(n)})=0,\label{eq30}\end{equation} с формальным решением $\;v={\Phi},$ $\;{\varphi}=x^{\rho_0}\,{\Phi},$ порядок решения $\;\Phi\,$ равен нулю. При этом, после степенного преобразования \eqref{eq29} каждое слагаемое уравнения \eqref{eq25} вида
$$\psi(x)\,x^{q_1}u^{q_{20}}(xu')^{q_{21}}\dots(x^nu^{(n)})^{q_{2n}}$$ переходит в сумму слагаемых вида
 \begin{equation}\label{eq33}
 \psi(x)\,x^{q_1+q_2\rho_0}v^{\tilde{q}_{20}}(xv')^{\tilde{q}_{21}}\dots (x^nv^{(n)})^{\tilde{q}_{2n}},
\end{equation}
где $ \tilde{q}_{20},\dots,\tilde{q}_{2n}\in\mathbb{Z}_+,\;$ $\tilde{q}_{20}+{\dots}+\tilde{q}_{2n}=q_2.\;$
 Следовательно, под действием степенного преобразования \eqref{eq29} вся  укороченная сумма $\;\hat{g}(x,u,u'\dots,u^{(n)})\;$ переходит в выражение $\;x^c\;\mathcal{P}_0(x,v,v'\dots,v^{(n)}),\;$ где $c$ -- это минимум всех значений $q_1+q_2\rho_0$, $\;\mathcal{P}_0(x,v,v'\dots,v^{(n)})$ -- это многочлен переменных $v,\,v'\,\dots,\,v^{(n)}\,$  с  коэффициентами в виде функций или формальных рядов нулевого порядка, а остальные слагаемые уравнения \eqref{eq25} перейдут в сумму слагаемых вида \,\eqref{eq33}\, с \,$\;q_1+q_2\rho_0>c.\;$ Поэтому уравнение \eqref{eq30} имеет вид
$$
  x^{c}\;[\mathcal{P}_0(x,v_0,\dots,v_n)+x^{\nu_1}\mathcal{P}_1(x,v_0,\dots,v_n)+{\dots}+x^{\nu_t}\mathcal{P}_t(x,v_0,\dots,v_n)]=0,
$$
где $v_j=x^jv^{(j)}$, $\mathcal{P}_0(x,v_0,\dots,v_n)$, $\dots,$ $\mathcal{P}_t(x,v_0,\dots,v_n)$ -- это многочлены переменных $v_0,\dots,v_n$ с коэффициентами в виде функций или формальных рядов нулевого порядка, числа $\nu_1,$ $\dots,$ $\nu_n\in\mathbb{R}_+,$ $\,\nu_1,$ $\dots,$ $\nu_n\neq 0.$ Разделим это уравнение  на $\;x^c,\;$ получим уравнение
\begin{equation}\label{eq46}
  \mathcal{P}_0(x,v_0,\dots,v_n)+x^{\nu_1}\mathcal{P}_1(x,v_0,\dots,v_n)+{\dots}+x^{\nu_t}\mathcal{P}_t(x,v_0,\dots,v_n)=0,
\end{equation} которое обозначим $\widetilde{G}(x,v_0,\dots,v_n)=0.$
Из записи \eqref{eq46} и предложений 1 -- 3 работы [\ref{OnClasses}] \, следует, что порядок $\;p\left(\widetilde{G}(x,\check{\Phi}_0,\dots,\check{\Phi}_n)\right)\geqslant 0\;$ (это следует из того, что при $\;\mathcal{P}_j\left(x,\check{\Phi}_0,\dots,\check{\Phi}_n\right)\neq 0\;$  порядок ряда $\;p\left(\mathcal{P}_j(x,\check{\Phi}_0,\dots,\check{\Phi}_n)\right)=0$, поскольку каждое слагаемое ряда $\;\mathcal{P}_j\left(x,\check{\Phi}_0,\dots,\check{\Phi}_n\right)\;$ является функцией нулевого порядка).

 Формальное решение $$\qquad u=\check{\Phi},\qquad\check{\Phi}=\sum\limits_{k_0=0}^{\infty}\tilde{\Phi}_{k_0}(x_1,\dots,x_N)\,\Phi_{k_0}(x)$$ уравнения \eqref{eq46}  запишем в виде
\begin{equation}\label{eq48}
\qquad u=\hat{\Phi}+w,\qquad\hat{\Phi}=\tilde{\Phi}_0(x_1,\dots,x_N)\,\Phi_0(x).
\end{equation}
Подставим  решение \eqref{eq48} в уравнение \eqref{eq46}, а затем (учитывая то, что уравнение \eqref{eq46}
 полиномиально по $v_0,\dots,v_n$) в выражении после подстановки раскроем скобки по формуле бинома Ньютона, получим выражение
$$\mathcal{P}_0(x,\hat{\Phi}_0,\dots,\hat{\Phi}_n)+\left[\frac{\partial\mathcal{P}_0(x,\hat{\Phi}_0,\dots,\hat{\Phi}_n)}{\partial v_0}\,w_0+{\dots}+\frac{\partial\mathcal{P}_0(x,\hat{\Phi}_0,\dots,\hat{\Phi}_n)}{\partial v_n}\,w_n\right]+$$
\begin{equation}
+{\dots}+x^{\nu_1}\mathcal{P}_1(x,\hat{\Phi}_0,\dots,\hat{\Phi}_n)+\dots=0,\qquad \displaystyle \hat{\Phi}_j=x^j\,\hat{\Phi}^{(j)},\qquad w_j=x^j\,w^{(j)}.
\label{eq47}
\end{equation}

Отдельно оценивая порядок каждого слагаемого уравнения \eqref{eq47},\, убеждаемся,\, что наименьший  порядок функций в нуле и в бесконечности имеют члены ряда
  $\;\mathcal{P}_0(x,\hat{\Phi}_0,$ $\dots,$ $\hat{\Phi}_n)$, при этом порядок ряда $p(\mathcal{P}_0(x,\hat{\Phi}_0,$ $\dots,$ $\hat{\Phi}_n))$
   либо равен нулю, либо равен $-\infty$ (если $\mathcal{P}_0(x,\hat{\Phi}_0,$ $\dots,$ $\hat{\Phi}_n)=0$).
    А порядки остальных ненулевых слагаемых равенства \eqref{eq47} больше нуля. Действительно, порядок
    $\;p(w_j)>p(\hat{\Phi}_j)=0$,\,  формальные ряды \,$\mathcal{P}_0(x,\hat{\Phi}_0,\dots,\hat{\Phi}_n),$\, $\dots,$\,
    $\mathcal{P}_t(x,\hat{\Phi}_0,\dots,\hat{\Phi}_n)$ и \,их  производные\, по\, $v_0,$\, $\dots,$\, $v_n$
    либо равны нулю, либо имеют порядок, равный нулю \, в силу предложений\, 2\, и\, 3 из [\ref{OnClasses}],\, \,   а числа $\nu_n>\dots>\nu_1>0$.
 Но поскольку ряд $\check{\Phi}$ формально удовлетворяет уравнению  \eqref{eq46}, то сумма всех слагаемых минимального порядка в выражении \eqref{eq47}
  должна быть равна нулю, т.~е.  $\mathcal{P}_0(x,\hat{\Phi}_0,\dots,\hat{\Phi}_n)=0$.
Делая в уравнении $\,\mathcal{P}_0(x, v_0,\dots,v_n)=0\,$ обратное преобразование к преобразованию \eqref{eq29} и домножая
результат на $x^{\,c},$ получаем, что
$$\hat{g}(x,\hat{\varphi}_0,\dots,\hat{\varphi}_n)=0,\qquad \hat{\varphi}_j=x^j\,\hat{\varphi}^{(j)},\qquad \hat{\varphi}=x^{\rho_0}\,\hat{\Phi}.\quad \Box$$

\bigskip

\bigskip
{\bf 3. Переход к уравнению специального вида.} Приведем уравнение \eqref{eq1} к полиномиальному виду. Для этого домножим обе части уравнения на множитель $x^2(x-1)^2y(y-1)(y-x)$,\, и перенесем правую часть уравнения в левую. Получим алгебраическое обыкновенное дифференциальное уравнение
\begin{equation}
2x^2(x-1)^2y(y-1)(y-x)y''-x^2(x-1)^2(3y^2-2xy-2y+x)y'^2+\label{eq3},
\end{equation}
$$+2xy(x-1)(y-1)(2xy-x^2-y)y'-2y^6a+4a(x+1)y^5-$$
$$-2\left((a+d)x^2+(4a+b+c-d)x+(a-c)\right)y^4+$$
$$+
4x\left((a+b+c+d)x+(a+b-c-d)\right)y^3-$$
$$-2\left((b+c)x^3+(a+4b-c+d)x^2+(b-d)x\right)y^2+4bx^2(x+1)y-2bx^3=0,$$
которое имеет те же самые решения, что и уравнение \eqref{eq1}, см. [\ref{MMO}].

Из работы [\ref{MMO}] следует, что каждое первое слагаемое $\,\varphi_0(x)\,$ экзотических или сложных формальных решений \eqref{eq2} шестого уравнения Пенлеве -- это рациональная функция переменной $\,\xi$,\, где $\,\xi=\ln^{-1}{x}\,$ либо $\,\xi=x^{\,{\rm i}\,\theta}$.\, Всякая рациональная функция от таких $\xi$ является функцией нулевого порядка по переменной $x$ в открытых секторах комплексной плоскости с вершинами в нуле и в бесконечности, которые не содержат в себе нули и полюсы этой функции. Позже мы, конечно, выпишем точные формулы для первых слагаемых  формальных решений вида \eqref{eq2} шестого уравнения Пенлеве.

Сделаем в уравнении \eqref{eq3} преобразование
\begin{equation}
y=\varphi_0(x)+x u.\label{eq4}
\end{equation}
После чего получаем \, {\it уравнение специального вида}
\begin{equation}\label{eq5}
L\left(\Phi_0,\dot{\Phi}_0,\ddot{\Phi}_0, U\right)+x M\left(x,\Phi_0,\dot{\Phi}_0,\ddot{\Phi}_0, U\right)+H\left(x,\Phi_0,\dot{\Phi}_0,\ddot{\Phi}_0\right)=0,
\end{equation}
где  $\;U=(u_0,u_1,u_2),$ $\;u_j=x^j u^{(j)},\,$ $$\,\Phi_0=\varphi_0(x),\;\dot{\Phi}_0=x\frac{d\varphi_0(x)}{dx},\;\ddot{\Phi}_0=x^2\frac{d^2\varphi_0(x)}{dx^2},$$  дифференциальный полином
$$
  L\left(\Phi_0,\dot{\Phi}_0,\ddot{\Phi}_0, U\right)=2\Phi_0^2\left(\Phi_0-1\right)u_2\,+\,
$$
\begin{equation}2\Phi_0\left(3\Phi_0^2-3\Phi_0\dot{\Phi}_0-3\Phi_0+2\dot{\Phi}_0\right)u_1\,-\,\label{eq6}\end{equation}
$$2\left(6a\Phi_0^5-10a\Phi_0^4+4a\Phi_0^3
  -4c\Phi_0^3-\Phi_0^3-3\Phi_0^2\ddot{\Phi}_0+3\Phi_0\dot{\Phi}_0^2+\Phi_0^2+2\Phi_0\ddot{\Phi}_0-\dot{\Phi}_0^2\right)u_0,$$
   $\;M\left(x,\Phi_0,\dot{\Phi}_0,\ddot{\Phi}_0, U\right)\,$ и $\,H\left(x,\Phi_0,\dot{\Phi}_0,\ddot{\Phi}_0\right)\,$ также являются дифференциальными полиномами. Уравнение \eqref{eq5} имеет решение
     \begin{equation}\label{eq7}
       u=\sum\limits_{k=0}^{\infty} \varphi_{k+1}(x)\,x^{k}.
     \end{equation}

{\bf 4. Анализ сложных разложений.}
 Далее будем исследовать уравнение \eqref{eq5} c формальным решением вида \eqref{eq7}, пользуясь результатами теории аналитических дифференциальных уравнений.
 Линейная часть уравнения \eqref{eq5} по $u_0,$ $u_1,$ $u_2$ состоит из членов, входящих в $L\left(\Phi_0,\dot{\Phi}_0,\ddot{\Phi}_0,\right.$ $\left. U\right)\,$
 и $\,xM(x,\Phi_0,\dot{\Phi}_0,\ddot{\Phi}_0,$ $U)$.
 При этом ведущими членами в окрестности точки $x=0$ являются члены, образующие
  $L\left(\Phi_0,\dot{\Phi}_0,\ddot{\Phi}_0, U\right)$. Это следует из того, что в $L$ порядки функций, стоящих перед $u_0,\,u_1$ и $u_2$ равны нулю,
  а в $M$ -- больше нуля.

 Согласно работе [\ref{MMO}] при $ a\neq c$, $\;a,\,c\neq 0$ шестое уравнение Пенлеве имеет формальное решение вида \eqref{eq2}
 с  коэффициентами $\varphi_k(x)$, которые являются рядами Лорана с конечными главными частями от $\ln^{-1}x$ и
 начальным членом
\begin{equation}\label{eq11}
  \varphi_0(x)=\frac{2(c-a)}{(c-a)^2(\ln x+C)^2-2a}.
\end{equation}

При подстановке формального решения \eqref{eq7} с асимптотикой \eqref{eq11}  в уравнение \eqref{eq5} получается формальный ряд
$$\sum\limits_{k=1}^{\infty}\left(L_k(\Phi_0,\dot{\Phi}_0,\ddot{\Phi}_0,\Phi_k,\dot{\Phi}_k,\ddot{\Phi}_k)-N_k(\Phi_0,\dot{\Phi}_0,\dots,\ddot{\Phi}_0,\Phi_{k-1},\dot{\Phi}_{k-1},\ddot{\Phi}_{k-1})\right)x^k,$$
 где \begin{equation}\label{eq17}
\Phi_j=\varphi_j(x),\;\displaystyle \dot{\Phi}_j=x\frac{d\varphi_j(x)}{dx},\;\displaystyle \ddot{\Phi}_j=x^2\frac{d^2\varphi_j(x)}{dx^2},
\end{equation}
\begin{equation}L_k\left(\Phi_0,\dot{\Phi}_0,\ddot{\Phi}_0,\Phi_k,\dot{\Phi}_k,\ddot{\Phi}_k\right)=\label{eq35}\end{equation}
$$=L\left(\Phi_0,\dot{\Phi}_0,\ddot{\Phi}_0,\Phi_k,\dot{\Phi}_k+(k-1)\Phi_k,\ddot{\Phi}_k+2(k-1)\dot{\Phi}_k+(k-1)(k-2)\Phi_k\right),$$
 функции  $N_k\left(\Phi_0,\dot{\Phi}_0,\ddot{\Phi}_0,\right.$ $\dots,$ $\Phi_{k-1},$ $\dot{\Phi}_{k-1},$ $\left.\ddot{\Phi}_{k-1}\right)$
 являются полиномами.

Таким образом, получаем, что уравнения на коэффициенты $\varphi_k(x)$ с натуральным $k$ имеют вид
  \begin{equation}\label{eq13}
  L_k\left(\Phi_0,\dot{\Phi}_0,\ddot{\Phi}_0,\Phi_k,\dot{\Phi}_k,\ddot{\Phi}_k\right)=
  N_k\left(\Phi_0,\dot{\Phi}_0,\ddot{\Phi}_0,\dots,\Phi_{k-1},\dot{\Phi}_{k-1},\ddot{\Phi}_{k-1}\right).
  \end{equation}

 Сделаем в уравнениях \eqref{eq13} переход от переменной $x$ к переменной $\chi$, где \begin{equation}\chi=(\xi+C)^{-1}=(\ln x+C)^{-1},\quad C\in\mathbb{C}.\label{eq7n}\end{equation}
    При этом отметим, что  преобразование \eqref{eq7n} в алгебраической дифференциальной сумме  не выводит ее из своего класса,
    поскольку $$x\,\frac{d y}{dx}= -\chi^2\frac{dy}{d\chi},\qquad
   x^2\,\frac{d^2 y}{dx^2}= \chi^4\frac{d^2y}{d\chi^2}+(2\chi^3+\chi^2)\frac{dy}{d\chi}.$$

Получившиеся уравнения домножим на $\displaystyle \frac{(-2a\chi^2+a^2-2ac+c^2)^4}{(-c+a)^{6}\chi^4}$. После этого они принимают вид
\begin{equation}\label{eq10}
\mathcal{L}_k\left(\chi,\hat{\varphi}_k(\chi),\frac{d\hat{\varphi}_k(\chi)}{d\chi},\frac{d^2\hat{\varphi}_k(\chi)}{d\chi^2}\right)=\mathcal{N}_k(\chi),
\end{equation}
где $$\mathcal{L}_k\left(\chi,\hat{\varphi}_k(\chi),\frac{d\hat{\varphi}_k(\chi)}{d\chi},\frac{d^2\hat{\varphi}_k(\chi)}{d\chi^2}\right)=$$
$$\chi^4 P_2(\chi)\frac{d^2\hat{\varphi}_k(\chi)}{d\chi^2}+\chi^2 P_1(\chi)\frac{d\hat{\varphi}_k(\chi)}{d\chi}+P_0(\chi)\hat{\varphi}_k(\chi),$$
$\hat{\varphi}_k(\chi)=\varphi_k(x),$
$P_2$ и $P_0$ -- многочлены четвертой степени, $P_1$ -- многочлен пятой степени, $P_2(0),$ $P_1(0),$ $P_0(0)\neq 0,$
а $\mathcal{N}_k(\chi)$ -- формальный ряд Лорана с конечной главной частью по $\chi$.

Заметим, что линейные дифференциальные операторы этих уравнений  в точке $\chi=0$ не являются фуксовыми. Это означает, что каждый из них имеет вид
  $$a_2(\chi)\chi^2\frac{d^2}{d\chi^2}+a_{1}(\chi)\chi\frac{d}{d\chi}+a_0(\chi),$$ где $a_2(\chi),$ $a_{1}(\chi),$  $a_0(\chi)$ -- голоморфные функции,
   и порядок нуля функции $a_2(\chi)$ строго меньше порядка нуля либо функции $a_1(\chi)$, либо функции $a_0(\chi)$.
Следовательно, формальные ряды ${\varphi}_k(x)$, которые являются рядами Лорана по $\ln^{-1}(x)$ c конечными главными частями, могут расходиться.
Тот факт, что каждый коэффициент является рядом Лорана доказывается просто  подстановкой формального ряда \eqref{eq7} с неопределенными коэффициентами
в уравнение \eqref{eq13} (где $\varphi_0(x)$ из \eqref{eq11}). Поскольку ведущее слагаемое в левой части уравнения \eqref{eq10} алгебраическое,
т.е. $$\mathcal{L}_k\left(\chi,\hat{\varphi}_k(\chi),\frac{d\hat{\varphi}_k(\chi)}{d\chi},\frac{d^2\hat{\varphi}_k(\chi)}{d\chi^2}\right)=-8 k^2 \hat{\varphi}_k(\chi)+\dots,$$
а другая часть уравнения -- это полином от предыдущих коэффициентов (формальных рядов Лорана от $\chi$),
то каждое такое уравнение имеет формальное решение в виде ряда Лорана от $\chi$. Поэтому $\varphi_k(x)$ -- суть формальные ряды Лорана от $\ln^{-1}x$.

При $\;a=c\neq 0\;$ шестое уравнение Пенлеве имеет также два семейства сложных формальных решений вида \eqref{eq2} c асимптотикой
\begin{equation}\label{eq100}
\varphi_0(x)=\frac{1}{\sqrt{2a}\;(\ln x+C)},\;\;C\in\mathbb{C}.
\end{equation}
 Исследования дальнейших коэффициентов этих формальных решений приводят
 к аналогичным результатам, что и для семейства формальных решений вида \eqref{eq2} c асимптотикой \eqref{eq11}. Соответствующие рассуждения здесь не приводим.

 Резюмируя результаты этого пункта получаем нижеследующее утверждение.

 \bigskip
 {\bf Утверждение 1.} {\it Коэффициенты $\,\varphi_k(x)\,$  $\,(k\geqslant 1)\,$ формальных решений  \eqref{eq2} c асимптотиками \eqref{eq11} и \eqref{eq100}
 шестого уравнения Пенлеве являются формальными рядами Лорана  с конечными главными частями от $\,\ln^{-1}x$,\, которые могут расходиться.}

\bigskip
Методы вычисления сложных формальных решений алгебраических обыкновенных дифференциальных уравнений предложены в [\ref{PreprBr}] и [\ref{ArticlBr}].

 Данная работа является частью будущей диссертации автора и поэтому исследования проводились только в рамках задач,
 включенных в диссертацию. Но отметим, что имеется множество других задач, которые возникают в этой теме.  В частности,
  интересно получить оценки порядков Жевре (оценки скорости роста коэффициентов)  рядов Лорана $\varphi_k(x)$ от $\ln^{-1}x$ и,
 если возможно, доказать их сходимость или расходимость и найти сумму ряда.

\bigskip

{\bf 5. Анализ экзотических разложений.} Далее мы будем исследовать формальные решения вида \eqref{eq7} уравнения \eqref{eq5} с коэффициентами
$\varphi_k(x)$ в виде рядов Лорана по $x^{{\rm i}\theta}$ с конечными главными частями
  и формальной асимптотикой
 \begin{equation}\label{eq15}
  \frac{-4C_0(2a-2c+C_1)}{x^{\sqrt{2c-2a-C_1}}(C_1^2+8C_1a+16a^2-16ac)-2C_1C_0+C_0^2x^{-\sqrt{2c-2a-C_1}}},
 \end{equation}
 где $\,C_0,\,C_1$ -- комплексные произвольные постоянные, $C_0\neq 0,\,$ $2c-2a-C_1\in\mathbb{R},$ $\,2c-2a-C_1<0,$ $\,\tau={\rm sgn}({\rm Im}\,\sqrt{2c-2a-C_1})$.

 При этом из предыдущих работ автора (см., например, [\ref{MMO}]) известны пять пар семейств базовых экзотических
 формальных решений $\mathcal{B}_0^{\tau}$, $\mathcal{B}_1^{\tau}$,  $\mathcal{B}_2^{\tau}$, $\mathcal{B}_6^{\tau}$ и $\mathcal{B}_7^{\tau}$.
  Оказалось, что все они являются подсемействами двупараметрического (по $C_0$ и $C_1$) семейства экзотических формальных решений
  вида \eqref{eq2} с асимптотикой \eqref{eq15} уравнения \eqref{eq1}.
  А именно, получаем семейства
\begin{itemize}
  \item $\,\mathcal{B}_0^{\tau}\,$ с асимптотикой
  \begin{equation}\label{eq15_B0}
 \varphi_0(x)=\frac{2c-C_3}{2a}\frac{1}{\alpha\cos^2(\ln(C_2x)^{\sqrt{C_3-2c}/2})+\beta\sin^2(\ln(C_2x)^{\sqrt{C_3-2c}/2})},
 \end{equation}
     при $a\neq 0,$ $\displaystyle C_0= \frac{\sqrt{C_3^2+4C_3a+4a^2-16ac}}{C_2^{\sqrt{2c-C_3}}}$, $\;C_1=C_3-2a$, где
   $\;C_3^2+4C_3a+4a^2-16ac\neq 0,$ $\;C_3\neq 2c,\;$
   $\;C_2\neq 0$, $\;2c-C_3\in\mathbb{R},\;$ $2c-C_3<0$,
  $\alpha$ и $\beta$ -- корни квадратного уравнения $2a t^2+(C_3-2a)t+2c-C_3=0$, $\;\tau=$ ${\rm sgn}({\rm Im}\sqrt{2c-C_3})$;

\item $\,\mathcal{B}_1^{\tau}\,$ с асимптотикой
 \begin{equation}\label{eq15_B1}
 \varphi_0(x)=\frac{1-\sqrt{c/a}}{1-C_2x^{\sqrt{2c}-\sqrt{2a}}},
 \end{equation}
  при $a\neq c\neq 0,$ $\displaystyle C_0=8\sqrt{a}(\sqrt{c}-\sqrt{a})\,C_2$, $C_1=4\sqrt{a}(\sqrt{c}-\sqrt{a})$, где
  $C_2\neq 0$, ${\rm Re}(\sqrt{2c}-\sqrt{2a})=0$, $\tau={\rm sgn}({\rm Im}(\sqrt{2c}-\sqrt{2a}))$;

\item $\,\mathcal{B}_2^{\tau}\,$ с асимптотикой
 \begin{equation}\label{eq15_B2}
 \varphi_0(x)=\frac{1+\sqrt{c/a}}{1-C_2x^{\sqrt{2c}+\sqrt{2a}}},
 \end{equation}
при $a\neq c\neq 0,$ $\displaystyle C_0=-8\sqrt{a}(\sqrt{c}+\sqrt{a})\,C_2$, $C_1=-4\sqrt{a}(\sqrt{c}+\sqrt{a})$, где $C_2\neq 0$,
${\rm Re}(\sqrt{2c}+\sqrt{2a})=0$, $\tau={\rm sgn}({\rm Im}(\sqrt{2c}+\sqrt{2a}))$;

 \item $\,\mathcal{B}_6^{\tau}\,$ с асимптотикой
 \begin{equation}\label{eq15_B6}
 \varphi_0(x)=\frac{1}{1+C_2 x^{\sqrt{2a}}},
 \end{equation}
 при $a\neq0,$ $c=0$, $C_0=8aC_2,$ $C_1=-4a$, где $C_2\neq 0$,   $\tau={\rm sgn}({\rm Im}\sqrt{2a})$;

 \item $\,\mathcal{B}_7^{\tau}\,$ с асимптотикой
 \begin{equation}\label{eq15_B7}
  \varphi_0(x)=\frac{2c-C_1}{C_1}\frac{1}{\sin^2(\ln(C_2x)^{\sqrt{C_1-2c}/2})}
 \end{equation}
 при $a=0,$ $\displaystyle C_0=-C_1/C_2^{\sqrt{2c-C_1}}$, $\;C_2\neq 0$, $\;2c-C_1\in\mathbb{R},\;$ $2c-C_1<0$,
  $\;\tau={\rm sgn}({\rm Im}\sqrt{2c-C_1})$.
   \end{itemize}

 Далее мы будем рассматривать самое общее двупараметрическое (по $C_0$ и $C_1$) семейство  экзотических формальных решений уравнения \eqref{eq1},
 которые имеют вид \eqref{eq2} с асимптотикой \eqref{eq15}.
 Отметим, что функция \eqref{eq15} является рациональной относительно $Cx^{{\rm i}\theta}$.

Сначала будем исследовать коэффициент $\varphi_1(x)$. Для этого изучим уравнение, которому он удовлетворяет. Это уравнение получается после подстановки
 экзотического формального решения \eqref{eq7} в уравнение \eqref{eq5}, упорядочивании членов согласно росту их порядков, выписывании членов первого порядка и деления их на $x$.
 Оно имеет вид
\begin{equation}{L}_1\left(\Phi_0,\dot{\Phi}_0,\ddot{\Phi}_0,\Phi_1,\dot{\Phi}_1,\ddot{\Phi}_1\right)+N_1\left(\Phi_0,\dot{\Phi}_0,\ddot{\Phi}_0\right)=0,\label{eq14}
  \end{equation}
где $\Phi_0,\dot{\Phi}_0,\ddot{\Phi}_0,\Phi_1,\dot{\Phi}_1,\ddot{\Phi}_1$ определены формулой \eqref{eq17},
$$L_1\left(\Phi_0,\dot{\Phi}_0,\ddot{\Phi}_0,\Phi_1,\dot{\Phi}_1,\ddot{\Phi}_1\right)=2\Phi^2_0(\Phi_0-1)\ddot{\Phi}_1+
2\Phi_0(3\Phi_0^2-3\Phi_0\dot{\Phi}_0-3\Phi_0+2\dot{\Phi}_0)\dot{\Phi}_1$$ $$-2(6a\Phi_0^5-10a\Phi_0^4+4a\Phi_0^3-4c\Phi_0^3+\Phi_0^3-3\Phi_0^2\ddot{\Phi}_0
+3\Phi_0\dot{\Phi}_0^2
+\Phi_0^2+2\Phi_0\ddot{\Phi}_0-\dot{\Phi}_0^2)\Phi_1,$$
$$N_1\left(\Phi_0,\dot{\Phi}_0,\ddot{\Phi}_0\right)=4a\Phi_0^5-2(4a+b+c-d)\Phi_0^4+4(a+b-c-d)\Phi_0^3-6\Phi_0^3\dot{\Phi}_0$$
$$-4\Phi_0^3\ddot{\Phi}_0+6\Phi_0^2\dot{\Phi}_0^2-2(b-d)\Phi_0^2+6\Phi_0^2\dot{\Phi}_0+2\Phi_0^2\ddot{\Phi}_0-2\Phi_0\dot{\Phi}_0^2+2\Phi_0\ddot{\Phi}_0-\dot{\Phi}_0^2.$$

Сделаем в выражении \eqref{eq15} и  уравнении \eqref{eq14}  переход от переменной $x$ к переменной $\chi=C\xi=Cx^{{\rm i}\theta}$,
 $\;\theta\in\mathbb{R}$, $\;\theta\neq 0$. При этом, отметим, что
$$x\,\frac{d y}{dx}= {\rm i}\theta\chi\frac{dy}{d\chi},\qquad
   x^2\,\frac{d^2 y}{dx^2}= -\chi^2\theta^2\frac{d^2y}{d\chi^2}-\theta(\theta+{\rm i})\chi\frac{dy}{d\chi},$$
т.~е. под действием данного преобразования  алгебраическая дифференциальная сумма переходит в алгебраическую дифференциальную сумму.
Здесь также $\chi$ используется вместо $\xi$ только для того, чтобы формулы были менее громоздкими.

Функция \eqref{eq15}, а следовательно и функции $\Phi_0,$ $\dot{\Phi}_0,$ $\ddot{\Phi}_0$ являются рациональными по $Cx^{{\rm i}\theta}$ и
в новой переменной $\chi=Cx^{{\rm i}\theta}$ имеют вид
$$\Phi_0=\frac{4\theta^2\chi}{A\chi^2+B\chi+1},$$ $$A=\theta^4+4(a+c)\theta^2+4(a-c)^2, B=2 \theta^2-4 (a- c),$$
\begin{equation}\label{eq39}
\dot{\Phi}_0 = -\frac{4\,{\rm i}\, \theta^3 \chi (A\chi^2-1)}{(A\chi^2+B\chi+1)^2},
\end{equation}
$$\ddot{\Phi}_0 = \frac{4\theta^3\chi(A^2({\rm i}-\theta)\chi^4+AB({\rm i}+\theta)\chi^3+6A\theta\chi^2-B({\rm i}-\theta)\chi-{\rm i}-\theta)}{(A\chi^2+B\chi+1)^3}.$$

Уравнение \eqref{eq14} для коэффициента $y=\varphi_1(x)$ в новой переменной и после умножения на выражение
$\displaystyle -\frac{(A\chi^2+B\chi+1)^6}{16\,\theta^4\,\chi^2}$ принимает вид
\begin{equation}\label{eq18}
\sum\limits_{j=0}^8\left(p_{2j}\, \chi^{j+2}\,\frac{d^2y}{d\chi^2}+p_{1j}\,\chi^{j+1}\,\frac{dy}{d\chi}+p_{0j}\,\chi^j y\right)+t_{j}\,\chi^{j}=0,
\end{equation}
где числа $p_{2j},$ $p_{1j},$ $p_{0j}, \;t_j\in\mathbb{C},$ $p_{20}$ и $p_{28}\neq 0.$
Поскольку уравнение \eqref{eq18} линейное неоднородное, то его решение имеет вид $$y=C_1\, y_1(\chi)+C_2\, y_2(\chi)+y_3(\chi),$$ $C_1,\,C_2$ --
произвольные комплексные постоянные, $y_1(\chi),\;y_2(\chi)$, $\;y_3(\chi)$ -- базисные решения.

Отметим, что наименьшее общие кратное уравнения \eqref{eq14} после замены переменной $\chi=Cx^{{\rm i}\theta}$ -- это многочлен $(A\chi^2+B\chi+1)^6$,
 при умножении на который,
 а затем делении результата на $-16\,\theta^4\,\chi^2$  получаем уравнение \eqref{eq18},
при
старшей производной $\;\displaystyle \frac{d^2y}{d\chi^2}\;$ которого стоит выражение
\begin{equation}\label{eq101}
\displaystyle -\frac{(A\chi^2+B\chi+1)^6}{8\theta^4}\Phi_0^2(\Phi_0-1).
\end{equation}
 Выражение  \eqref{eq101} -- это многочлен десятой степени по $\chi$, который имеет 5 корней $\chi=0,\;a_1,\;a_2,\;a_3,\;a_4$ без учета их кратности. А именно,
 $\chi=0$ -- корень кратности 2, $\chi=a_1$ и $\chi=a_2$ -- корни кратности 1, где
 $a_1$ и $a_2$ -- это решения уравнения $(A-4\theta^2)\chi^2+B\chi+1=0,$ и $a_3$ и $a_4$ -- корни кратности 3, где $a_3$ и $a_4$ --
 это решения уравнения $A\chi^2+B\chi+1=0.$
  Поэтому уравнение \eqref{eq18} имеет шесть особых точек $\chi=0,\infty,$ $a_1,$ $a_2,$ $a_3,$ $a_4.$
Так как коэффициенты $p_{20}$ и $p_{28}\neq 0$, то линейный дифференциальный оператор уравнения \eqref{eq18}  фуксов в нуле, и в бесконечности.
  Сделав в уравнении \eqref{eq18} замены переменных $\chi=\zeta+a_j$,\, нетрудно проверить, что и в этих точках
  линейный дифференциальный оператор уравнения \eqref{eq18}  является фуксовым. Локальный вид решений фуксовых уравнений,
  которые являются линейными однородными обыкновенными дифференциальными уравнениями, известен см. [\ref{Bolibruch}].
 Что касается линейного неоднородного обыкновенного дифференциального уравнения с фуксовой однородной частью,
 то мы воспользовались методом вариаций произвольных постоянных, чтобы получить локальный вид его частного решения.
 Из этого всего следует, что любое решение уравнения \eqref{eq18} локально в окрестности каждой особой точки имеет вид
$$y=\sum\limits_{i=1,2}C_iF_i(\chi-a_j)(\chi-a_j)^{\lambda_i}\,\ln^{\mu_i}{(\chi-a_j)}+$$
\begin{equation}\label{eq19}
+F_3(\chi-a_j)(\chi-a_j)^{\lambda_3}\,\ln^{\mu_3}{(\chi-a_j)},
\end{equation}
где $F_1(\chi), F_2(\chi), F_3(\chi)\in\mathbb{C}\{\chi\},$ $\lambda_1,\lambda_2,\lambda_3\in\mathbb{C},$ $\mu_1,\;\mu_2,\;\mu_3\in\mathbb{Z}.$

  Напомним, что каждое формальное решение линейного обыкновенного
  дифференциального уравнения с фуксовым оператором сходится. Кроме того, особенности решений линейного обыкновенного
  дифференциального уравнения могут быть только в особенностях уравнения [\ref{Golubev}].

  Далее мы найдем все формальные решения уравнения \eqref{eq18} во всех  его особых точках.
  Исследование формальных решений в окрестностях особых точек уравнения
  позволит нам определить, что $y_3(\chi)$ -- рациональная функция, а $y_1(\chi)$ и
  $y_2(\chi)$ -- неоднозначные функции, которые в точках $\chi=0$ и $\chi=\infty$ имеют трансцендентное ветвление.

  \bigskip
  {\bf Утверждение 2.} {\it Коэффициент $\varphi_1(x)$ формального решения \eqref{eq7} уравнения \eqref{eq5} с асимптотикой \eqref{eq15}
   является рациональной функцией от $Cx^{{\rm i}\theta}$. А именно, $\varphi_1(x)=y_3(Cx^{{\rm i}\theta})$, где $y_3(\chi)$ --
   это частное рациональное решение уравнения \eqref{eq18}}.

\bigskip
{\bf Доказательство утверждения 2.}
   При исследовании решений уравнения \eqref{eq18} мы будем пользоваться тем, что все разложения, соответствующие решениям вида \eqref{eq19},
    содержатся в классе формальных решений, которые могут быть найдены с помощью методов плоской степенной геометрии.
   Согласно теореме 1, если некоторое такое разложение является решением уравнения \eqref{eq18},
   то его укорочение обязательно является решением укороченного уравнения.
   Укороченные уравнения соответствуют ребрам или вершинам многоугольника Ньютона-Брюно.
   Поэтому рассмотрев все решения укороченных уравнений, соответствующих всем ребрам и вершинам (которых всегда конечное число),
   выбрав среди них те, которые являются первыми слагаемыми формальных решений (формальными асимптотиками),
   и затем вычислив характеристические числа (из-за которых могут появляться в разложении нецелые степени и логарифмы)
   дифференциального оператора ведущих членов линейной части уравнения,  мы определим  характер особенностей решения уравнения \eqref{eq18}.

  Полностью выписывать уравнение \eqref{eq18} мы не будем, поскольку оно черезвычайно громоздкое.
  Многоугольник Ньютона-Брюно уравнения \eqref{eq18} -- это прямоугольник с вершинами $[(0,0),\;(0,1),\;(8,1),\;(8,0)]$.
  Согласно этому многоугольнику вблизи нуля на роль асимптотик нам подходят решения укороченных уравнений с порядком по переменной $\chi$
  меньшим нуля (если уравнение соответствует верхней левой вершине)  и равным нулю (если уравнение соответствует левому вертикальному ребру),
  а вблизи бесконечности -- с порядком по переменной $\chi$ большим нуля (если уравнение соответствует правой верхней вершине)
  и равным нулю (если уравнение соответствует правому вертикальному ребру).

  Выпишем укороченные (приближенные) уравнения при $\chi\rightarrow 0$ и при $\chi \rightarrow\infty$. При $\chi\rightarrow 0$ имеем два
  линейных уравнения с постоянными коэффициентами
  $$\left(-2\theta^2\chi^2\frac{d^2y}{d\chi^2}+2\theta(\theta+2{\rm i})\chi\frac{dy}{d\chi}-2(\theta +{\rm i})^2y\right)=0$$ и
$$-\left(2\theta^2\chi^2\frac{d^2y}{d\chi^2}+2\theta(\theta+2{\rm i})\chi\frac{dy}{d\chi}-2(\theta +{\rm i})^2y+(\theta+{\rm i})^2-1 +2b-2d\right)=0,$$
где левая часть первого уравнения -- это однородная часть левой части второго уравнения.
При этом, однородное уравнение соответствует вершине $(0,1)$, а неоднородное -- ребру $[(0,1),\;(0,0)]$. Решение неоднородного уравнения имеет вид
  \begin{equation}\label{eq20}
 y=(C_1+C_2\ln{\chi})\chi^{\scriptstyle 1+\frac{\scriptstyle \rm i}{\scriptstyle \theta}}+\frac{(\theta+{\rm i})^2+1+2b-2d}{2(\theta+{\rm i})^2},
\end{equation}
а однородного --
$$y=(C_1+C_2\ln{\chi})\chi^{\scriptstyle 1+\frac{\scriptstyle \rm i}{\scriptstyle \theta}},$$ где $C_1$, $C_2$ -- произвольные комплексные постоянные.
Решение однородного уравнения не подходит на роль формальной асимптотики, поскольку имеет порядок по переменной $\chi$  больший нуля,
а решение неоднородного уравнения \eqref{eq20} годится на роль формальной асимптотики, поскольку имеет порядок по переменной $\chi$ равный нулю.

Далее пользуясь методами степенной геометрии [\ref{BrunoUMN}], находим, что
формальная асимптотика \eqref{eq20} продолжается в формальный ряд
\begin{equation}\label{eq23}
y=\left(C_1\chi\sum\limits_{k=0}^\infty a_{1k} \chi^k+C_2\chi\ln{\chi} \sum\limits_{k=0}^\infty a_{2k} \chi^k\right)\chi^{\frac{\scriptstyle \rm i}{\scriptstyle \theta}}+\sum\limits_{k=0}^\infty a_{3k} \chi^k,
\end{equation}
$a_{1k},\,a_{2k},\,a_{3k}\in\mathbb{C},$ $\,a_{10}=a_{20}=1$, $$a_{30}=\frac{(\theta+{\rm i})^2+1+2b-2d}{2(\theta+{\rm i})^2},$$
удовлетворяющий уравнению \eqref{eq18}.

При $\chi\rightarrow\infty$ также имеем два линейных уравнения с постоянными коэффициентами
$$A^4\chi^8\left(-2\theta^2\chi^2\frac{d^2y}{d\chi^2}-2\theta(3\theta -2{\rm i})\chi\frac{d y}{d\chi}-2(\theta-{\rm i})^2y\right)=0$$
и
$$A^4\chi^8\left(-2\theta^2\chi^2\frac{d^2y}{d\chi^2}-2\theta(3\theta -2{\rm i})\chi\frac{d y}{d\chi}-2(\theta-{\rm i})^2y-(\theta-{\rm i})^2-1+2b-2d\right)=0,$$
где левая часть первого уравнения -- это однородная часть левой части второго уравнения.
    Однородное уравнение соответствует вершине $(8,1)$, а неоднородное -- ребру $[(8,1),\;(8,0)]$. Решение неоднородного уравнения имеет вид
  \begin{equation}\label{eq22}
 y=(C_1+C_2\ln{\chi})\chi^{\scriptstyle -1+\frac{\scriptstyle \rm i}{\scriptstyle \theta}}+\frac{(\theta-{\rm i})^2+1+2b-2d}{2(\theta-{\rm i})^2},
\end{equation}
а однородного --
$$y=(C_1+C_2\ln{\chi})\chi^{\scriptstyle -1+\frac{\scriptstyle \rm i}{\scriptstyle \theta}},$$ где $C_1$, $C_2$ -- произвольные комплексные постоянные.
Решение однородного уравнения не подходит на роль формальной асимптотики, поскольку имеет порядок по переменной $\chi$ меньший нуля,
а решение неоднородного уравнения \eqref{eq22} годится на роль формальной асимптотики, поскольку имеет порядок по переменной $\chi$ равный нулю.

Формальная асимптотика \eqref{eq22} продолжается в формальный ряд
\begin{equation}\label{eq24}
y=\left(\frac{C_1}{\chi}\sum\limits_{k=0}^\infty \frac{b_{1k}}{\chi^k}+\frac{C_2}{\chi}\ln{\chi} \sum\limits_{k=0}^\infty \frac{b_{2k}} {\chi^k}\right)\chi^{\frac{\scriptstyle \rm i}{\scriptstyle \theta}}+\sum\limits_{k=0}^\infty \frac{b_{3k}}{\chi^k},
\end{equation}
$b_{1k},\,b_{2k},\,b_{3k}\in\mathbb{C},$ $\,b_{10}=b_{20}=1$, $$b_{30}=\frac{(\theta-{\rm i})^2+1+2b-2d}{2(\theta-{\rm i})^2},$$ удовлетворяющий уравнению \eqref{eq18}.

Теперь рассмотрим уравнение \eqref{eq18} в окрестностях точек $a_1$ и $a_2$, которые являются корнями уравнения $(A-4\theta^2)\chi^2+B\chi+1=0$.
Для этого в уравнении \eqref{eq18} сделаем замены независимой переменной $\chi=\zeta+a_j$, $j=1,\;2$.
Для каждого корня получается свое линейное неоднородное уравнение с фуксовой однородной частью.
Поскольку для каждого из этих уравнений дальнейшие рассуждения и вычисления аналогичны, мы будем писать о них, как об одном уравнении.
Получившееся после замены переменной уравнение можно записать в виде
\begin{equation}\label{eq18_1}
\sum\limits_{j=0}^8\left(P_{2j}\, \zeta^{j+1}\,\frac{d^2y}{d\zeta^2}+P_{1j}\,\zeta^{j}\,\frac{dy}{d\zeta}+P_{0j}\,\zeta^j y+T_{j}\,\zeta^{j}\right)=0,
\end{equation}
где $P_{2j},$ $P_{1j},$ $P_{0j},\;T_j\in\mathbb{C},$ $P_{20}\neq 0.$
Многоугольник Ньютона-Брюно уравнения \eqref{eq18_1} -- это четырехугольник с вершинами $(-1,1),\;(0,0),\;(8,0),\;(8,1)$.
Нас интересуют формальные решения, соответствующие укороченным уравнениям, решения которых являются кандидатами в асимптотики в окрестности точки $\zeta=0$,
т.~е. формальные решения, соответствующие левой вершине $(-1,1)$ с порядком меньшим единицы и левому ребру $[(-1,1),\;(0,0)]$ с порядком равным единице. Решая соответственные укороченные уравнения
 $-\zeta y''+y'=0$ и $-\zeta y''+y'=\alpha$, $\;\alpha\in\mathbb{C}$, $\;\alpha$
 рационально зависит от коэффициентов уравнения \eqref{eq18_1}, получаем асимптотики
 $\,y=C_1\neq 0\,$ и $\,y=\alpha\zeta$, $\,C_1\,$ -- произвольная комплексная постоянная.
 Этим двум асимптотикам соответствует одно двупараметрическое по $\,C_1,\,C_2\,$ ($C_2\in\mathbb{C}$)\, семейство формальных решений
\begin{equation}\label{eq31}
y=C_1\sum\limits_{k=0}^\infty c_{1k} \zeta^k+C_2\zeta^2\sum\limits_{k=0}^\infty c_{2k} \zeta^k+\zeta\sum\limits_{k=0}^\infty c_{3k} \zeta^k,
\end{equation}
$c_{1k},\,c_{2k},\,c_{3k}\in\mathbb{C},$ $\,c_{10}=c_{20}=1$, $c_{30}=\alpha,$  уравнения \eqref{eq18_1}.

\bigskip
Теперь осталось исследовать уравнение \eqref{eq18} в окрестностях остальных двух особых точек $a_3$ и $a_4$, которые являются корнями уравнения
 $\,A\chi^2+B\chi+1=0$. Для этого в уравнении \eqref{eq18} сделаем замены независимой переменной $\chi=\zeta+a_j$, $j=3,\;4$.
 И снова для каждого корня получается свое линейное неоднородное уравнение с фуксовой однородной частью.
  Но для каждого из этих уравнений дальнейшие рассуждения и вычисления аналогичны, поэтому мы будем писать о них, как об одном уравнении.
  Получившееся после замены переменной уравнение можно записать в виде
\begin{equation}\label{eq18_1_n}
\sum\limits_{j=0}^8\left(S_{2j}\, \zeta^{j+3}\,\frac{d^2y}{d\zeta^2}+S_{1j}\,\zeta^{j+2}\,\frac{dy}{d\zeta}+S_{0j}\,\zeta^{j+2} y+K_{j}\,\zeta^{j}\right)=0,
\end{equation}
где $S_{2j},$ $S_{1j},$ $S_{0j},\;K_j\in\mathbb{C},$ $S_{20}\neq 0.$
 Многоугольник Ньютона-Брюно этого уравнения -- это четырехугольник с вершинами $(1,1),\;(0,0),\;(8,0),\;(8,1)$.
Нас интересуют формальные решения, соответствующие укороченным уравнениям, решения которых являются кандидатами в асимптотики в окрестности точки $\zeta=0$,
т.~е. формальные решения, соответствующие левой вершине $(1,1)$ с порядком меньшим минус единицы и левому ребру $[(1,1),\;(0,0)]$
с порядком равным минус единице. Решая соответственные укороченные уравнения
 $\zeta^2(\zeta y''+3y')=0$ и $\zeta^2(\zeta y''+3y')=\alpha$, $\;\alpha\in\mathbb{C}$, $\;\alpha$  рационально зависит от коэффициентов уравнения \eqref{eq18_1},
  получаем асимптотики $\displaystyle \,y=\frac{C_1}{\zeta^2},$ $C_1\neq 0\,$ и $\displaystyle\,y=\frac{\alpha}{\zeta}$,
  $\,C_1\,$ -- произвольная комплексная постоянная.
  Этим двум асимптотикам соответствует одно двупараметрическое по $\,C_1,\,C_2\,$ ($C_2\in\mathbb{C}$)\, семейство формальных решений
\begin{equation}\label{eq32}
y=\frac{C_1}{\zeta^2}\sum\limits_{k=0}^\infty d_{1k} \zeta^k+C_2\sum\limits_{k=0}^\infty d_{2k} \zeta^k+\frac{1}{\zeta}\sum\limits_{k=0}^\infty d_{3k} \zeta^k,
\end{equation}
$d_{1k},\,d_{2k},\,d_{3k}\in\mathbb{C},$ $\,d_{10}=d_{20}=1$, $d_{30}=\alpha,$  уравнения \eqref{eq18_1}.

Еще раз отметим, что в окрестности всех особых точек $0,\,\infty,\,a_1,\,a_2,\,a_3,\,a_4$ уравнения \eqref{eq18} его формальные решения сходятся.
Анализируя вид этих локальных решений можно утверждать, что общее решение $y=C_1y_1(\chi)+C_2y_2(\chi)+y_3(\chi)$  этого уравнения имеет вид
$$
y=C_1f_1(\chi)\,\chi^{\frac{\scriptstyle \rm i}{\scriptstyle \theta}}+C_2\ln{\chi}\, f_2(\chi)\,\chi^{\frac{\scriptstyle \rm i}{\scriptstyle \theta}}+f_3(\chi),
$$
где $\;f_1(\chi),$ $f_2(\chi)$, $f_3(\chi)\;$ -- рациональные функции от $\chi$, т.~е. функции $\;y_1(\chi)=C_1f_1(\chi)\,\chi^{\frac{\scriptstyle \rm i}{\scriptstyle \theta}}\;$
 и $\;y_2(\chi)=C_2\ln{\chi}\, f_2(\chi)\,\chi^{\frac{\scriptstyle \rm i}{\scriptstyle \theta}}\;$ -- это
неоднозначные функции, которые имеют трансцендентное ветвление в нуле и в бесконечности, и поэтому в этих точках не могут представляться локально рядами Лорана или Тейлора.
А функция $y_3(\chi)=f_3(\chi)$ -- рациональная функция, которая во всех особых точках уравнения представляется степенными рядами с целыми показателями степени.
Из этого следует, что коэффициент $\varphi_1(x)=y_3(Cx^{{\rm i}\theta}).$ $\Box$

\bigskip
Далее будем изучать свойства всех коэффициентов ${\varphi}_k(x)$ формального решения \eqref{eq7} уравнения \eqref{eq5} с асимптотикой \eqref{eq15}.
При подстановке этого формального решения в уравнение получается формальный ряд
$$\sum\limits_{k=1}^{\infty}\left(L_k(\Phi_0,\dot{\Phi}_0,\ddot{\Phi}_0,\Phi_k,\dot{\Phi}_k,\ddot{\Phi}_k)-N_k(\Phi_0,\dot{\Phi}_0,\dots,\ddot{\Phi}_0,\Phi_{k-1},\dot{\Phi}_{k-1},\ddot{\Phi}_{k-1})\right)x^k,$$
 где
$\Phi_j,$ $\dot{\Phi}_j$, $\ddot{\Phi}_j$ определены формулой \eqref{eq17},
$L_k\left(\Phi_0,\dot{\Phi}_0,\ddot{\Phi}_0,\Phi_k,\dot{\Phi}_k,\ddot{\Phi}_k\right)$ определены формулой \eqref{eq35},
 $N_k\left(\Phi_0,\dot{\Phi}_0,\ddot{\Phi}_0,\right.$ $\dots,$ $\Phi_{k-1},$ $\dot{\Phi}_{k-1},$ $\left.\ddot{\Phi}_{k-1}\right)$
 являются полиномами своих аргументов.

Из этого уравнения следует, что если ряд \eqref{eq7} удовлетворяет уравнению \eqref{eq5}, то каждый коэффициент ${\varphi}_k(x)$,  $k\in\mathbb{N}$
удовлетворяет уравнению
\begin{equation}L_k(\Phi_0,\dot{\Phi}_0,\ddot{\Phi}_0,\Phi_k,\dot{\Phi}_k,\ddot{\Phi}_k)=
N_k(\Phi_0,\dot{\Phi}_0,\dots,\ddot{\Phi}_0,\Phi_{k-1},\dot{\Phi}_{k-1},\ddot{\Phi}_{k-1}).\label{eq36}\end{equation}
Обращаем внимание читателя, что явное выражение левой части уравнения \eqref{eq36} нам известно для любого $k$,
 а явно выписывать правую часть можно лишь последовательно, вычисляя предыдущие коэффициенты.

Перейдем  от переменной $x$ к переменной $\chi=C\xi=Cx^{{\rm i}\theta}$,
 $\;\theta\in\mathbb{R}$, $\;\theta\neq 0$, в коэффициентах $\varphi_k(x)=\hat{\varphi}_k(\chi)$, $k\in\mathbb{N}$, и соответствующих им уравнениях \eqref{eq36}.
 После чего, уравнения \eqref{eq36}
 принимают вид
 \begin{equation}\label{eq34}
 Q_{2}(\chi)\;\chi^2\,\frac{d^2\hat{\varphi}_k(\chi)}{d\chi^2}
+Q_{1k}(\chi)\;\chi\,\frac{d\hat{\varphi}_k(\chi)}{d\chi}+Q_{0k}(\chi)\hat{\varphi}_k(\chi)=\mathcal{N}_k(\chi),
\end{equation}
где
$$Q_2(\chi)=-2\theta^2\,\Phi_0^2(\Phi_0-1),$$
$$Q_{1k}(\chi)=2\,{\rm i}\,\theta\,({\rm i}\,\theta+2k-3)Q_2(\chi)+{\rm i}\,\theta\,Q_1(\chi),$$
$$Q_{0k}(\chi)=(k-2)(k-1)Q_2(\chi)+(k-1)Q_1(\chi)+Q_0(\chi),$$
$$Q_{1}(\chi)=2\Phi_0(3\Phi_0^2-3\Phi_0\dot{\Phi}_0-3\Phi_0+2\dot{\Phi}_0),$$
$$Q_{0}(\chi)=-2(6a\Phi_0^5-10a\Phi_0^4+4a\Phi_0^3-4c\Phi_0^3+\Phi_0^3-3\Phi_0^2\ddot{\Phi}_0
+3\Phi_0\dot{\Phi}_0^2
+\Phi_0^2+2\Phi_0\ddot{\Phi}_0-\dot{\Phi}_0^2)\Phi_1,$$
$\Phi_0,\;\dot{\Phi}_0,\;\ddot{\Phi}_0$ определены формулой \eqref{eq39},
а $\mathcal{N}_k(\chi)$ -- формальный ряд Лорана с конечной главной частью по $\chi$.

\bigskip
{\bf Утверждение 3.} {\it
Уравнение \eqref{eq34} имеет решение $\hat{\varphi}_k(\chi)\,$  в виде формального ряда Лорана с конечной главной частью (или формального ряда Тейлора) по $\chi$.
 }

\bigskip
{\bf Доказательство утверждения 3.} Рассмотрим для каждого значения $k$ уравнение \eqref{eq34}. Будем искать решения этих уравнений в виде рядов
 \begin{equation}\label{eq41}
 \hat{\varphi}_k(\chi)=\chi^{r_k}\sum\limits_{j=0}^{\infty}\alpha_{kj}\chi^j,
 \end{equation}
 $r_k\in\mathbb{Z}$, $\alpha_{kj}\in\mathbb{C}$.
  Случай $k=1$ был рассмотрен в утверждении 1.
Рассмотрим случай $k=2$. В этом случае правая часть уравнения \eqref{eq34} -- это полином от $\chi$  и
 рациональных по  $\chi$ функций $\Phi_0$ и $\Phi_1$, т.~е.
 в проколотой окрестности нуля $N_2$ представляется сходящимся рядом Лорана
  $\chi^{R_2}\sum\limits_{j=0}^{\infty}A_{2j}\chi^j$, 
 $R_2\in\mathbb{Z}$, $A_{2j}\in\mathbb{C}$.
 Разложим в ряды Лорана коэффициенты $Q_2,\;Q_{1k},\;Q_{2k}$, которые тоже являются рациональными функциями по $\chi$,
 и подставим ряд
$\displaystyle
 \hat{\varphi}_2(\chi)=\chi^{r_2}\sum\limits_{j=0}^{\infty}\alpha_{2j}\chi^j
$
  в уравнение \eqref{eq34},
 получим выражение
 \begin{equation}\label{eq40_2}
     \sum\limits_{j=0}^{\infty}\left(\frac{8}{B^2a^2}(2{\rm i}-\theta)^2\theta^4\chi^2+O(\chi^3)\right)\alpha_{2j}\chi^{r_2+j}=
     \chi^{R_2}\sum\limits_{j=0}^{\infty}A_{2j}\chi^j.
 \end{equation}
Здесь $\;a\neq 0,\; B\neq 0$, $\;2{\rm i}-\theta\neq 0$, так как $\theta\in\mathbb{R}$.
Всегда можно так подобрать числа $r_2, R_2$, что ряд в левой и правой частях равенства \eqref{eq40_2} начинается с одной и той же степени. Тогда
каждый коэффициент $\alpha_{2j}$  однозначно определяется.   Пусть теперь $k=3$, тогда $N_3$ -- это полином от $\chi$, (рациональных по $\chi$ функций)\,$\Phi_0$ и  $\Phi_1$,  и (ряда Лорана) $\Phi_2$.
При подстановке ряда $\eqref{eq41}$ в уравнение \eqref{eq34}  получаем уравнение
     \begin{equation}\label{eq40}
     \sum\limits_{j=0}^{\infty}\left(\frac{8}{B^2a^2}(k{\rm i}-k\theta+\theta)^2\theta^4\chi^2+O(\chi^3)\right)\alpha_{kj}\chi^{r_k+j}=
     \chi^{R_k}\sum\limits_{j=0}^{\infty}A_{kj}\chi^j,
     \end{equation}
где $R_k\in\mathbb{Z}$, $A_{kj}\in\mathbb{C}$. Рассуждения в этом случае аналогичны рассуждениям в случае $k=2$. Из
них следует, что  $\hat{\varphi}_3(\chi)$ также является рядом Лорана с конечной главной частью. Продолжая по индукции доказательство,
приходим к выводу, что каждое уравнение \eqref{eq34} имеет решение  в виде ряда Лорана с конечной главной частью по $\chi$.
$\Box$

\bigskip
{\bf Утверждение 4.} {\it Каждый коэффициент $\varphi_k(x)=\hat{\varphi}_k(\chi)\,$ разложения \eqref{eq7} уравнения \eqref{eq5} с асимптотикой \eqref{eq15}
  в проколотой окрестности точки $\chi=0$ представляется сходящимся рядом Лорана с конечной главной частью (или рядом Тейлора) по $\chi$.}

  {\bf Доказательство утверждения 4} следует из предыдущего утверждения и из того,
  что линейный дифференциальный оператор каждого уравнения \eqref{eq34} фуксов в нуле. $\Box$


\bigskip
{\bf Утверждение 5.}
{\it Уравнение \eqref{eq34} имеет   шесть особых точек $0,$ $\infty,$ $a_1,$ $a_2,$ $a_3,$ $a_4\in\mathbb{C}$, в каждой из которых линейный
 дифференциальный оператор уравнения фуксов.}

{\bf Доказательство утверждения 5.} Левая часть уравнения имеет шесть особых точек $0,$ $\infty,$ $a_1,$ $a_2,$ $a_3,$ $a_4\in\mathbb{C}$.
В каждой из которых соответствующий ей линейный дифференциальный оператор фуксов.
 А правая часть уравнения \eqref{eq34} при $k=2$ это полином от рациональных функций с теми же самыми особыми точками.
 Следовательно, решение уравнения \eqref{eq34}
 при $k=2$ имеет те же самые особые точки, что и уравнение. Теперь если $k=3$, то правая часть уравнения \eqref{eq34} -- это полином
 от рациональных функций и функции (возможно с трансцендентным ветвлением), но с теми же самыми особыми точками. Продолжая наши рассуждения,
 убеждаемся, что это верно для любого номера $k$. $\Box$

\bigskip
Из утверждения 5 следует, что решения уравнений \eqref{eq34} могут иметь особенности только в точках
$\chi=0,$ $\infty,$ $a_1,$ $a_2,$ $a_3,$ $a_4\in\mathbb{C}$. Это означает, что если существует настоящее решение уравнения \eqref{eq5}
с заданным разложением  \eqref{eq7} с асимптотикой \eqref{eq15}, то особенности такого решения могут содержатся и накапливаться только в секторах с биссектрисами
$Cx^{{\rm i}\theta}=0,$ $\infty,$ $a_1,$ $a_2,$ $a_3,$ $a_4$.

Отметим, что радиус сходимости каждого ряда Лорана $\hat{\varphi}_k(\chi)$  параметрически зависит от положений
особых точек $\chi=a_1,$ $a_2,$ $a_3,$ $a_4\in\mathbb{C}$ уравнения \eqref{eq34}.

Автором этой работы также было доказано, что  коэффициенты $\varphi_3(x)$ и $\varphi_4(x)$ разложения \eqref{eq7} с асимптотикой \eqref{eq15} уравнения \eqref{eq5}
также являются рациональными функциями от $Cx^{{\rm i}\theta}$. Вычисления громоздки и не дают информации о дальнейших коэффициентах, поэтому здесь их не приводим.
Но все это позволяет выдвинуть гипотезу.

\bigskip
{\bf Гипотеза.} {\it Все коэффициенты $\varphi_k(x)=\hat{\varphi}_k(\chi)$ -- это рациональные функции от $\chi=Cx^{{\rm i}\,\theta}$.}

\bigskip\bigskip
\begin{center}
СПИСОК ЛИТЕРАТУРЫ\end{center}

\begin{enumerate}

\item \label{Gromak}  \emph{Gromak I.V.,  Laine I.,
Shimomura S. } Painlev\'e Differential Equations in the Comp\-lex
Plain. Berlin, New York:  Walter de Gruyter. 2002.

\item\label{MMO} \emph{Брюно А.Д., Горючкина И.В.  } Асимптотические разложения решений шестого уравнения Пенлеве // Труды ММО.  2010. Т. 71. С.~6--118.

\item\label{OnClasses} {\it Горючкина И.В. } Классы формальных решений конечного порядка обыкновенного дифференциального уравнения // 	
Чебышевский сборник. 2016. Т. 17. № 2(58). С. 64-87

\item\label{Guzzetti} {\it Guzzetti D. } Poles Distribution of PVI transcendents close to a critical point // Physica D.
2012. doi:10.1016/j.physd.2012.02.015.

\item\label{OnConv} {\it Gontsov, R.R., Goryuchkina, I.V. } On the convergence of generalized power series satisfying an
algebraic ODE. Asympt. Anal. 2015. 93(4). P. 311--325.

\item\label{GG} {\it Gontsov R., Goryuchkina I. } An analytic proof of the Malgrange-Sibuya theorem
on the convergence of formal solutions of an ODE. J. Dynam. Control Syst. 2016. V. 22(1). P. 91-100.

\item\label{BrunoUMN}  \emph{Брюно А.Д. }  Асимптотики и разложения
решений обыкновенного дифференциального
    уравнения //  УМН. 2004. Т.~59. №~3. С.~31--80.

\item\label{PreprBr} \emph{Брюно А.Д. } О сложных разложениях
решений ОДУ // Препринты ИПМ им. М.В.Келдыша. 2011. №~15. 26 с.

\item\label{ArticlBr} {\it Брюно А.Д. } Элементы нелинейного анализа, Математический форум, серия ``Итоги науки. Юг России'', ЮМИВНЦ. 2015. C. 13--33.

\item\label{Bolibruch} {\it Болибрух А.А. } Обратные задачи монодромии в аналитической теории дифференциальных уравнениях. М.:МЦНМО. 2009. 200 с.

\item\label{Golubev} {\it Голубев В.В. } Лекции по аналитической теории дифференциальных уравнений.  М.-Л.: Гостехтеориздат. 1941. 400 с.

\end{enumerate}

\end{document}